\documentclass[12pt]{amsart}
\usepackage{amscd,amssymb}
\usepackage[arrow,matrix]{xy}

\theoremstyle{plain}
\newtheorem{thm}[subsection]{Theorem}
\newtheorem{lem}[subsection]{Lemma}
\newtheorem{prop}[subsection]{Proposition}
\newtheorem{cor}[subsection]{Corollary}

\theoremstyle{definition}
\newtheorem{rk}[subsection]{Remark}

\newtheorem{ex}[subsection]{Example}

\numberwithin{equation}{section}
\newcommand{\F}{{\mathcal F}}

\newcommand{\LL}{{\mathcal L}}

\newcommand{\V}{{\mathcal V}}

\newcommand{\Z}{\mathbb{Z}}

\newcommand{\C}{\mathbb{C}}

\newcommand{\PP}{\mathbb{P}}

\newcommand{\T}{\mathbb{T}}

\DeclareMathOperator{\Hom}{Hom}

\DeclareMathOperator{\im}{im}

\begin{document}

\title [On the irreducible components of characteristic varieties]
{On the irreducible components of characteristic varieties  }

\author[Alexandru Dimca]{Alexandru Dimca }
\address{  Laboratoire J.A. Dieudonn\'e, UMR du CNRS 6621,
                 Universit\'e de Nice-Sophia Antipolis,
                 Parc Valrose,
                 06108 Nice Cedex 02,
                 FRANCE.}
\email {dimca@math.unice.fr}



\begin{abstract}
This is a quick survey on the characteristic varieties associated to rank one local systems on a 
 smooth, irreducible, quasi-projective complex variety $M$. A key new result is Proposition \ref{prop1},
 giving additional information on the constructible sheaf $\F=R^0f_*(\LL)$, where $\LL$ is a rank one local system on $M$ and $f:M \to S$ is a surjective morphism with $S$ a smooth curve and the generic fiber $F$ of $f$ connected. Corollary \ref{cor2} says that for $S$ compact, the singular support of $\F$
 cannot be a singleton.
\end{abstract}

\maketitle

\section{} \label{s1}
Let $M$ be a smooth, irreducible, quasi-projective complex variety
and let $\T(M)=\Hom(\pi_1(M),\C^*)$ be the character variety of 
$M$. This is an algebraic group whose identity irreducible component is an algebraic torus
$(\C^*)^{b_1(M)}$. Clearly the group $\T(M)$ is connected if and only if the integral homology group $H_1(M)$ is torsion free.

 The {\em characteristic varieties}\/ of $M$ 
are the jumping loci for the cohomology of $M$, with 
coefficients in rank~$1$ local systems:
\begin{equation} 
\label{eq:charvarx}
\V^i_k(M)=\{\rho \in \T(M) \mid \dim H^i(M, \LL_{\rho})\ge k\}.
\end{equation}
When $i=1$, we use the simpler notation $\V_k(M)=\V^1_k(M)$.

Foundational results on the structure of the cohomology 
support loci for local systems on quasi-projective algebraic 
varieties were obtained by  Beauville \cite{Beau}, 
Green and Lazarsfeld \cite{GL}, Simpson \cite{Sim} (for the proper case), 
and Arapura \cite{A} (for the quasi-projective case and first characteristic varieties $\V_1(M)$). 

\begin{thm} \label{thm1}
The strictly positive dimensional irreducible components of the first characteristic variety $\V_1(M)$
are translated subtori in $\T(M)$ by elements of finite order. When $M$ is proper, then all the 
components of $\V_k^i(M)$ are translated subtori in $\T(M)$ by elements of finite order.
\end{thm}
In the non-proper case, we get unitary characters  $\rho_j \in \T(M)$ as isolated points of the first characteristic varieties $\V_1(M)$), see Arapura \cite{A}. Recent results by N. Budur give a precise description in terms of complex tori and convex rational polytopes of some related characteristic varieties, see  \cite{B}, and imply in particular that these unitary characters $\rho_j$'s are actually torsion points in $\T(M)$.

\medskip

In the sequel we  concentrate ourselves on the strictly positive dimensional irreducible components of the first characteristic variety $\V_1(M)$. They have the following rather explicit description, given by
Arapura \cite{A}.

\begin{thm} \label{thm2}
Let $W$ be a $d$-dimensional irreducible component of the first characteristic variety $\V_1(M)$, with $d>0$.
Then there is a regular morphism $f:M \to S$ onto a smooth curve $S$ with $b_1(S)=d$ such that the generic fiber $F$ of $f$ is connected, and a torsion character $\rho \in \T(M)$ such that the composition 
$$\pi_1(F) \stackrel{i_{\sharp}} \longrightarrow       \pi_1(M) \stackrel{\rho} \longrightarrow  \C^*,$$
where $i:F \to M$ is the inclusion, is trivial and
$$W=\rho \cdot f^*(\T(S)).$$

\end{thm}

If we fix a regular mapping $f:M \to S$ as above, one may ask how many irreducible components $W=\rho \cdot f^*(\T(S))$
can be obtained by varying the torsion character $\rho$. The answer, due to Beauville in the case $M$ proper \cite{Beau}, and to the present author in the general case \cite{D3},  is given by the following.

\begin{thm} \label{thm3}
For a given regular mapping $f:M \to S$ as above, the associated irreducible components $W=\rho \cdot f^*(\T(S))$
of $\V_1(M)$ are parametrized by the Pontrjagin dual ${\widehat T(f)}$ of the finite abelian group
$$T(f)=\frac {\ker \{f_*:H_1(M) \to H_1(S)\}}{\im \{i_*:H_1(F) \to H_1(M)\}}$$
 if $\chi(S)<0$ and by the non-trivial elements of this Pontrjagin dual ${\widehat T(f)}$ if $\chi(S)=0.$
\end{thm}
Note that $b_1(S)>0$ implies $\chi(S)\le 0$, hence we have a precise answer in all the possible cases.

\begin{cor} \label{cor1}
Let $W=\rho \cdot f^*(\T(S))$ be a strictly positive dimensional irreducible component of the first characteristic variety $\V_1(M)$ which is translated, i.e. $1 \notin W$. Then the associated subtorus $W'=f^*(\T(S))$ is also an irreducible component of the first characteristic variety $\V_1(M)$, unless $\chi(S)=0$, i.e. we are in one of the following 2 cases:

\noindent (i) $\dim W=2$ and $S$ is a compact smooth curve of genus 1, or

\noindent (ii) $\dim W=1$ and $S=\C^*$. 

\end{cor}

One may ask what about the strictly positive dimensional irreducible components of the higher characteristic varities $\V_q(M)$.
The answer is given by the following result, see \cite{D3}, Cor. 4.7.
For a constructible sheaf $\F$ on $S$, $\Sigma(\F)$ denotes the singular support of $\F$, i.e. the minimal finite subset $\Sigma$ in $S$ such that $\F$ is a local system on $S \setminus \Sigma$, see \cite{D1} for more details.

\begin{thm} \label{thm4}

If $W_{f,\rho}=\rho \otimes f^*(\T(S))$ is a strictly positive dimensional irreducible component of $\V_1(M)$,
 then
$$\dim H^1(M, \LL_{\rho} \otimes f^{-1}\LL) \ge -\chi(S)+  | \Sigma (R^0f_*(\LL_{\rho}))|   $$
with equality for all but finitely many local systems $\LL \in \T(S)$. In particular, $W_{f,\rho}$ is an irreducible component of $\V_q(M)$,
for any $1 \le q \le q(f,\rho):=-\chi(S)+ | \Sigma (R^0f_*(\LL_{\rho}))|$. Conversely, any positive
dimensional irreducible component of $\V_q(M)$ for $q \ge 1$ is of this type.

\end{thm}
In other words, this theorem tells us which is the generic dimension of $H^1(M,\LL')$ when $\LL'$
belongs to an irreducible component $W=W_{f,\rho}$ of $\V_1(M)$. Note that if $1 \in W_{f,\rho}$ (i.e. the character $\rho$ is trivial), then this generic dimension is precisely $-\chi(S)>0$.

Some of the main properties of the constructible sheaf $\F=R^0f_*(\LL_{\rho})$ are established in  \cite{D3}. In particular one has the following two results.

\begin{lem} \label{lem3}
Let $\LL_{\rho}$ be a rank one local system on $M$, $F$ the generic fiber of $f:M \to S$ and set 
$\F=R^0f_*(\LL_{\rho})$. Then either

\medskip

\noindent (i) the restriction $\LL_{\rho}|F$ is trivial, $\F | (S \setminus \Sigma(\F))$ is a rank one local system
and $\F_s =0$ if and only if  $s \in \Sigma(\F)$, or

\medskip

\noindent (ii) the restriction $\LL_{\rho}|F$ is non-trivial and $\F=0$.

\end{lem}

In the sequel we assume that the restriction $\LL_{\rho}|F$ is trivial, the other case being trivial.
In this case there is an induced character $\tilde \rho:T(f) \to \C^*$.

\begin{thm} \label{thm4.5}

Let $f:M \to S$ be a surjective morphism, with connected generic fiber $F$, and let $\tilde \rho:T(f) \to \C^*$
be the character associated to $\rho$.
Then the singular support $\Sigma (\F)$  of the associated constructible sheaf $\F$  is empty if and only if the character $\tilde \rho$ is trivial.

\end{thm}

Here we continue the study of this sheaf $\F$ by the following key result.

\begin{prop} \label{prop1}
For a given regular mapping $f:M \to S$ as above, and for a non-trivial element $\tilde \rho$ in the Pontrjagin dual ${\widehat T(f)}$, one has a natural adjunction isomorphism
$$\F=Rj_*j^{-1}\F$$
where $j:S \setminus \Sigma(\F) \to S$ is the inclusion. In particular, the local system
$j^{-1}\F$ on $S \setminus \Sigma(\F)$ is non trivial.

\end{prop}

\proof
It is known that a point  $c \in S$ is in $ \Sigma(\F)$  if and only if for a small disc $D_c$ centered
at $c$, the restriction of the local system $\LL_{\rho}$ to the associated tube $T(F_c)=f^{-1}(D_c)$ about the fiber $F_c$
is non-trivial. Let $T(F_c)'=T(F_c) \setminus F_c$ and note that the inclusion $i: T(F_c)' \to T(F_c)$
induces an epimorphism at the level of fundamental groups.

Hence, if $c \in  \Sigma(\F)$, then $\LL_{\rho}|T(F_c)' $ is a non-trivial rank one local system. In particular
$$H^0(T(F_c)',\LL_{\rho})=0.$$
If we apply the Leray spectral sequence to the locally trivial fibration
$$F \to T(F_c)' \to D_c'$$
where $D_c'=D_c \setminus \{c\}$, we get
$$H^0(D_c',\F)=H^0(T(F_c)',\LL_{\rho})=0.$$
It follows that $\F|D_c'$ is a non-trivial rank one local system. Hence $H^0(D_c',\F)=H^1(D_c',\F)=0$,
which proves the isomorphism $\F=Rj_*j^{-1}\F$.

\endproof 

\begin{cor} \label{cor2}
With the above notation, if $S$ is a compact curve, then $|\Sigma(\F)| \ne 1.$

\end{cor}

\proof

Assume that $S$ is a compact curve and that  $\Sigma(\F)=\{c\}$. Then one has the following.

\noindent (i) The cycle associated to the boundary of the small disc $D_c$ in the proof above is a trivial
element in the integral homology group $H_1(S \setminus \Sigma(\F))$. 

\noindent (ii) The local system $\F|(S \setminus \Sigma(\F))$ corresponds to a (non-trivial) homomorphism
$$H_1(S \setminus \Sigma(\F)) \to \C^*.$$

It follows that $\F|D_c'$ is a trivial local system, in contradiction with the proof of Proposition \ref{prop1}.

\endproof 

\begin{rk} \label{rk} If we assume that $M$ is proper, Corollary \ref{cor2} follows from Beauville \cite{Beau}.
Indeed, if the multiple fibers of $f$ are $m_1F_1$,...,$m_sF_s$, it follows from Remark 1.7 in \cite{Beau}
that
$${\widehat T(f)}=\{(\hat k_1,...,\hat k_s)\in \oplus _{i=1,s}\Z/m_i\Z~~|~~ \sum_{i=1,s}k_i/m_i \in \Z\}.$$
This implies that for $M$ proper and $s \le 1$, one has ${\widehat T(f)}=1$.
Note that Example \ref{exMASTER} shows that one may have $|\Sigma(\F)| = 1$ when $S$ is non compact. More precisely,
in this case there is just one multiple fiber, hence $s=1$, but ${\widehat T(f)}=\Z/2\Z$. Therefore there are
subtle differences between the proper and the non-proper cases.
\end{rk}

When $M$ is a hypersurface complement in some complex projective space $\PP^n$, e.g. $M$ is a hyperplane arrangement complement, then the following additional features hold,  \cite{D1}.

\begin{thm} \label{thm5}
Let $M$ be a hypersurface complement in a complex projective space $\PP^n$, and let
$W$ be a $d$-dimensional irreducible component of the first characteristic variety $\V_1(M)$, with $d>0$. Then the associated curve $S$ as in Theorem \ref{thm2} is obtained from $\C$ by deleting $d$
points. If $1 \in W$, then the generic dimension of $H^1(M,\LL')$ when $\LL' \in W$ is $d-1$.
\end{thm}

It follows that the regular mapping $f:M \to S$ associated to such a component $W$ is a pencil of hypersurfaces in $\PP^n$, a point of view which has led to new results on the irreducible components
of the first characteristic variety $\V_1(M)$, for which we refer to \cite{FY}, \cite{PY}.

In general, it is quite difficult to compute the irreducible components of the characteristic varieties,
especially those not passing through the origin.
Here is a beautiful and instructive example.

\begin{ex} \label{exMASTER} This is a key example discovered by A. Suciu, see Example 4.1 in  \cite{S1}
and Example 10.6 in  \cite{S2}. Consider the line arrangement in $\PP^2$ given by the equation
$$xyz(x-y)(x-z)(y-z)(x-y-z)(x-y+z)=0.$$

We number the lines of the associated affine arrangement in $\C^2$ (obtained by setting $z=1$) as follows: $L_1: x=0$, $L_2:x-1=0$, $L_3: y=0$, $L_4: y-1=0$, $L_5:x-y-1=0$, $L_6:x-y=0$ and $L_7:x-y+1=0$, see the pictures in  Example 4.1 in  \cite{S1} and Example 10.6 in  \cite{S2}. We consider also the line at infinity
$L_8:z=0$.
As stated in  Example 4.1 in  \cite{S1}, there are

\medskip

\noindent (i) Seven local components: six of dimension 2, corresponding to the triple points, and one of dimension 3, for the quadruple point.

\medskip

\noindent (ii) Five components of dimension 2, passing through 1, coming from the following neighborly
partitions (of braid subarrangements): $(15|26|38)$, $(28|36|45)$, $(14|23|68)$, $(16|27|48)$ and $(18|37|46)$. For instance, the pencil corresponding to the first partition is given by
$P=L_1L_5=x(x-y-z)$ and $Q=L_2L_6=(x-z)(x-y)$. Note that $L_3L_8=yz=Q-P$, is a decomposable fiber in this pencil.

\medskip

\noindent (iii) Finally, there is a 1-dimensional component $W$ in $\V_1(M)$ with 
$$\rho_W=(1,-1,-1,1,1,-1,1,-1) \in \T(M) \subset (\C^*)^8$$
and $f_W:M \to \C^*$ given by
$$f_W(x:y:z)=\frac{x(y-z)(x-y-z)^2}{(x-z)y(x-y+z)^2}$$
or, in affine coordinates
$$f_W(x,y)=\frac{x(y-1)(x-y-1)^2}{(x-1)y(x-y+1)^2}.$$

Then $W \subset \V_1(M)$ and $W \cap  \V_2(M)$ consists of two characters, $\rho_W$ above and
$$\rho _W'=(-1,1,1,-1,1,-1,1,-1).$$
Note that this component $W$ is a translated component corresponding to case (ii) in Corollary \ref{cor1}.

More precisely, it can be shown that $T(f_W)=\Z/2\Z$ and that the two characters $\rho_W$ and  $\rho_W'$ induce both the only non-trivial element in
 the Pontrjagin dual ${\widehat T(f_W)}$.
Indeed, corresponding mapping $f_W:M \to \C^*$ has just one bifurcation point, namely $c=\{1\}$. 
The corresponding fiber $F_c$ is just $2(L\cap M)$ where $L:x+y-1=0$ is exactly the line from the $B_3$-arrangement that was
deleted in order to get Suciu's arrangement.

\end{ex}

\bigskip

There are number of interesting applications of the characteristic varieties $\V_q(M)$ to the study of the
fundamental groups of smooth algebraic varieties, see \cite{DPS}, \cite{DPS1}, \cite{DPS2}.

\end{document}